\documentclass{article}%
\usepackage{amsmath}
\usepackage{amsfonts}
\usepackage{amssymb}
\usepackage{graphicx}
\usepackage{subfigure}
\usepackage{epsfig}%
\usepackage{color}
\usepackage{multirow}
\usepackage{setspace}
\usepackage{hyperref}
\usepackage{enumitem}
\usepackage{hyperref}
\usepackage{authblk}
\newtheorem{theorem}{Theorem}

\newtheorem{proposition}[theorem]{Proposition}

\begin{document}
\doublespacing

\title{\textbf{On one-sample Bayesian tests for the mean}}   

\author[1]{Ibrahim Abdelrazeq}
\author[2]{Luai Al-Labadi}
\author[3] { Ayman Alzaatreh}
\affil[1]{Department of Mathematics and Computer Science, Rhodes College, Memphis, TN 38134, USA.
E-mail: Abdelrazeqi@Rhodes.edu}
\affil[2]{Department of Mathematical \& Computational Sciences, University of Toronto Mississauga, Ontario L5L 1C6, Canada.
E-mail:  lallabadi@sharjah.ac.ae}
\affil[3]{Department of Mathematics and Statistics, American University of Sharjah, Sharjah, UAE.
E-mail:  ayman.alzaatreh@aus.edu}

%




\date{}
\maketitle

\pagestyle {myheadings} \markboth {} {On one-sample Bayesian tests}

\begin{abstract}
This paper deals with a new Bayesian approach to the standard one-sample $z$- and $t$- tests. More specifically, let $x_1,\ldots,x_n$ be an independent random sample from a normal distribution with mean $\mu$ and variance $\sigma^2$. The goal is to test the null hypothesis $\mathcal{H}_0: \mu=\mu_1$ against all possible alternatives.   The approach is based on using the well-known formula of the Kullbak-Leibler divergence between two normal distributions (sampling and hypothesized distributions selected in an appropriate way). The change of the distance from a priori to a posteriori is compared through the relative belief ratio (a measure of evidence). Eliciting the prior, checking for prior-data conflict and bias are also considered. Many theoretical properties of the procedure have been developed. Besides it's simplicity, and unlike the classical approach, the new approach possesses attractive and distinctive features  such as giving evidence in favor of the null hypothesis.  It also avoids several undesirable paradoxes, such as Lindley's paradox that may be encountered by some existing Bayesian methods. The use of the approach has been illustrated through several examples. 

\par

 \vspace{9pt} \noindent\textsc{Keywords:}  Hypothesis testing, Kullbak-Leibler divergence, one-sample $t$-test, one-sample $z$-test, relative belief inferences.

 \vspace{9pt}

\noindent { \textbf{MSC 2000}} 62F03, 62F15

\end{abstract}
	
\section{Introduction}
The one-sample hypothesis testing  is a primary topic in any introductory statistics course. It involves  the selection of a reference value  $\mu_1$ for the (unknown) population mean $\mu$. More specifically,  let $x=(x_1,\ldots,x_n)$ be an independent random sample taken from $N(\mu, \sigma^2)$, where
$\sigma^2$ is the population variance. The interest is to test the hypothesis $H_0:\mu = \mu_1$, where $\mu_1$
is a given real number. Within the classical frequentist frame work, if $\sigma$ is known, then the $z$-test is commonly used for testing $H_0$ against the two-sided alternative $H_1: \mu \neq \mu_1$. The test statistics in this case is
$$z=\frac{\bar{x}-\mu_1}{\sigma/\sqrt{n}},$$
where $\bar{x}$ is the sample mean. For a significant level $\alpha$,  the critical value $Z_{\alpha/2}$ is defined to be the $1-\alpha/2$ quantile of the  standard normal distribution. Also, the $p$-value is equal to $2P(Z>|z|)$, where $Z$ has  the standard normal distribution. Then, $H_0$ is rejected if $|z| \ge Z_{\alpha/2}$ or the $p$-value less than $\alpha$. On the other hand, if $\sigma$ is unknown, then the test statistic is
$$t=\frac{\bar{x}-\mu_1}{s/\sqrt{n}},$$
where $s$ is the sample standard deviation.  For a test with significant level $\alpha$, let $t_{n-1,\alpha/2}$ be the $1-\alpha/2$ quantile of the $t$ distribution with $n-1$ degrees of freedom. The two sided $p$-value is equal to $2P(T>|t|)$, where $T$ has  the t-distribution with $n-1$ degrees of freedom. Similar to the $z$-test,  $H_0$ is rejected if $|t| \ge t_{n-1,\alpha/2}$ or the $p$-value is less than $\alpha$.

While the above approach for hypothesis testing is well-known and stable, it is difficult to find  an alternative Bayesian counterpart in the literature. An exception includes the work of Rouder, Speckman, Sun, and Morey (2009) who proposed a Bayesian test, where $\sigma$ is unknown, using the Bayes factor (ratio of the marginal densities of the two models; Kass and Raftery, 1995). They placed the Jeffreys prior for $\sigma$ and the Cauchy prior on $\mu/\sigma$.  They provided a web-based program (c.f. pcl.missouri.edu) in order to facilitate the use of their test. Remarkably, the authors mentioned detailed criticisms of using the $p$-values in hypothesis testing. For example, they indicated that the $p$-values do not allow researchers to state evidence for the null hypothesis. They also overstate the evidence against the null hypothesis. Although the $p$-value converges to zero as the sample size increases when the null hypothesis is false which is a desirable feature, the $p$-values are all equally likely and uniformly distributed between 0 and 1 when null is true. This distribution holds regardless of the sample size which means that increasing the sample size in this case will not help gaining evidence for the null hypothesis.  In fact, this reflects Fisher's sight that the null hypothesis can only be rejected and never accepted.  Other relevant work, but  in the two-sample problem set up, includes G\"onen,  Johnson, Lu and Westfall (2005) and Wang and Lui (2016). For more recent articles about the limitations of using $p$-values in hypotheses testing, we refer the reader to Evans (2015), Wasserstein and Lazar (2016), and  references therein.


Unlike the previous work, the hyperparameters of the prior in the new approached Bayesian are elicited and tested against prior-data conflict and against being biased.   For this, two elicitation algorithms developed by Evans (2015, 2018) are considered. In fact, the success of any Bayesian approach depends significantly on a proper selection of the hyperparameters of the prior. Part of the elicitation process involves checking the elicited prior for the prior-data conflict and the bias (see Section 2). Then the concentration of the distribution of the Kullbak-Leibler divergence between the prior and the model of interest is compared to that between the posterior and the model.  If the posterior is more concentrated about the hypothesized distribution than the prior, then this is evidence in favor of the null hypothesis and if the posterior is less concentrated then this is evidence against the null hypothesis. This comparison is made via a relative belief ratio, which measures the evidence in the observed data for or against the null. A measure of the strength of this evidence is also provided. So, the methodology is based on a direct measure of statistical evidence.  We point out that, relative belief ratios have been recently used in problems that involve goodness of fit test and model checking. See, for example,  Al-Labadi (2018), Al-Labadi and Evans (2018) and  Al-Labadi, Zeynep and Evans (2017, 2018) and Evans and Tomal (2018).

The proposed method brings many advantages  to the problem of hypothesis testing. Besides its simplicity, and unlike the classical approach, the new approach possesses attractive and desirable features such as giving evidence in favor of the null hypothesis. Also, checking the prior for bias and prior-data conflict permits avoid several undesirable paradoxes, such as Lindley's paradox that may be encountered by the standard Bayesian methods that are based, for instance, on the Bayes factor (Evans, 2015).

The remainder of this paper is organized as follows.  A general discussion about the relative belief ratio is given in Section 2. The definition and some fundamental properties of the Dirichlet process are presented in Section 3. In Section 4, an explicit expression to compute  Anderson-Darling distance between the Dirichlet process and its base measure is derived. In Section 5, a Bayesian nonparametric test  for assessing multivariate normality is discussed and some of its relevant properties are developed. A computational algorithm to calculate the relative belief ratio for the implementation of the proposed test is developed in Section 6. In Section 7, the performance of the proposed test is established via four simulated examples and two real data sets.  Finally, some concluding remarks are given in Section 8. All technical proofs are included in the supplementary material.

\section{Inferences Using Relative Belief}
Suppose we have a statistical model that is given by the density function $f_\theta(x)$ (with respect to some measure), where $\theta$ is an unknown parameter that belongs to the parameter space $\Theta$. Let $\pi(\theta)$ be the prior distribution of $\theta$. After observing the data $x$, by Bayes' theorem, the posterior distribution of $\theta$  is given by the density
\begin{equation*}
  \pi(\theta|x) = \frac{f_\theta(x)\pi(\theta)}{m(x)},
\end{equation*}
where
\begin{equation*}
m(x) = \int f_\theta(x)\pi(\theta) d\theta
\end{equation*}
is the prior predictive density  of the data.

Suppose that the interest  is to make inference about an arbitrary parameter $\psi=\Psi(\theta)$. Let $\Pi_{\Psi}$ denote the prior measure of $\psi$ with density $\pi_{\Psi}$. Let the corresponding posterior measure and density of $\psi$ be $\Pi_{\psi}(\cdot\,|\,x)$  and $\pi_{\Psi}%
(\cdot\,|\,x),$ respectively. The relative belief ratio for  a hypothesized value $\psi_0$ of $\psi$ is
 defined by $RB_{\Psi}(\psi_0\,|\,x)=\lim_{\delta\rightarrow0}\Pi_{\Psi
}(N_{\delta}(\psi_0\,)|\,x)/\Pi_{\Psi}(N_{\delta}(\psi_0\,))$, where $N_{\delta
}(\psi_0\,)$ is a sequence of neighbourhoods of $\psi$ converging   nicely (see, for example, Rudin (1974)) to
$\psi$ as $\delta\rightarrow0.$ When $\pi_{\Psi}$ and  $\pi_{\Psi}(\cdot\,|\,x)$ are continuous at $\psi,$
\begin{equation*}
RB_{\Psi}(\psi_0\,|\,x)=\pi_{\Psi}(\psi_0\,|\,x)/\pi_{\Psi}(\psi_0),
\end{equation*}
is the ratio of the posterior density to the prior density at $\psi_0.$ That is,
$RB_{\Psi}(\psi_0\,|\,x)$ is measuring how beliefs have changed that
$\psi_0$ is the true value from \textit{a priori} to \textit{a posteriori}. Baskurt and Evans (2013) proved that

\begin{equation}
RB_{\Psi}(\psi_0\,|\,x)=m_T(T(x)|\psi_0)/m_T(T(x)),
\label{savage1}%
\end{equation}
where  $T$ is a minimal sufficient statistic of the model and $m_T$ is
the prior predictive density of $T$.  The previous authors  referred to (\ref{savage1}) as the Savage-Dickey ratio. It is to be noted that a relative
belief ratio is similar to a Bayes factor (Kass and Raftery, 1995), as both are measures of evidence,
but the latter measures it via the change in an odds ratio.  A discussion about the relationship between relative belief ratios and Bayes factors is detailed in (Baskurt and Evans, 2013). More specifically, when a Bayes factor is defined via a limit in the continuous case, the limiting value is the corresponding relative belief ratio.

By a basic principle of evidence,  $RB_{\Psi}(\psi_0\,|\,x)>1$ means that the data led to an increase in the probability that $\psi_0$ is correct, and so there is
evidence in favour of $\psi_0,$ while $RB_{\Psi}(\psi_0\,|\,x)<1$ means that the data led
to a decrease in the probability that $\psi_0$ is correct, and so there is
evidence against $\psi_0$. Clearly, when $RB_{\Psi}(\psi_0\,|\,x)=1$, then there is no
evidence either way.


It is also important to calibrate whether this is strong
or weak evidence for or against $\mathcal{H}_{0}$. As suggested in Evans
(2015), a useful calibration of
$RB_{\Psi}(\psi_{0}\,|\,x)$ is obtained by computing the tail probability
\begin{equation}
\Pi_{\Psi}(RB_{\Psi}(\psi\,|\,x)\leq RB_{\Psi}(\psi_{0}\,|\,x)\,|\,x).
\label{strength}%
\end{equation}
One way to view (\ref{strength}) is as the posterior probability that the true value of $\psi$ has a relative
belief ratio no greater than that of the hypothesized value $\psi_{0}.$ When $RB_{\Psi}(\psi_{0}\,|\,x)<1,$ there is evidence
against $\psi_{0},$ then a small value for (\ref{strength}) indicates a large
posterior probability that the true value has a relative belief ratio greater
than $RB_{\Psi}(\psi_{0}\,|\,x)$ and so there is strong evidence against
$\psi_{0}.$ When $RB_{\Psi}(\psi_{0}\,|\,x)>1,$ there is evidence in favour
of $\psi_{0},$ then a large value for (\ref{strength}) indicates a small
posterior probability that the true value has a relative belief ratio greater
than $RB_{\Psi}(\psi_{0}\,|\,x))$. Therefore, there is strong evidence in favour of
$\psi_{0},$ while a small value of (\ref{strength}) only indicates weak
evidence in favour of $\psi_{0}.$

One of the key concerns with Bayesian inference methods is that the prior can bias the analysis. Following Evans (2015), let $M(\cdot|\psi)$ denote the conditional prior predictive distribution of the  data given that
$\Psi(\theta)=\psi$, so $M(A|\psi)=\int_{\Theta}\left(\int_{A}f_{\theta}(x)dx\right)\pi(\theta|\psi)d\theta$ is the conditional prior probability that the data is in the set $A$. The bias against $H_0: \Psi(\theta)=\psi_0$ can be measured by computing
\begin{equation}
M(RB_{\Psi}(\psi_0 | x) \le  1 | \psi_0)
\label{bias_against}%
\end{equation}
and this is the prior probability that evidence will be
obtained against $H_0$ when it is true. If the bias against $H_0$ is large, subsequently reporting, after seeing
the data, then there is evidence against $H_0$ is not convincing.On the other hand, the bias in favor of $H_0$ is given by
\begin{equation}
M(RB_{\Psi}(\psi_0 | x) \ge 1 | \psi^{\prime}_0)
\label{bias_favor}%
\end{equation}
for values $\psi_0 \neq \psi^{\prime}_0$ such that the difference between $\psi_0$ and $\psi^{\prime}_0$
 represents
the smallest difference of practical importance; note that this tends to decrease as $\psi^{\prime}_0$
 moves farther
away from $\psi_0$. When the bias in favor is large, subsequently reporting, after seeing the data, then the
is evidence in favor of $H_0$ is not convincing. 


Another concern regarding priors is to measure the compatibility  between the prior and the data. A chosen prior may be incorrect by being strongly contradicted by the data (Evans, 2015). A possible contradiction between the data and the prior is referred to as a prior-data conflict. If the prior primarily places its mass in a region of the parameter space where the data suggest the true value does not lie, then there is a \emph{prior-data conflict} (Evans and Moshonov, 2006). That is, prior-data conflict will occur whenever there is  only a tiny overlap between the effective support regions of the model and  the prior. In such situation, we must be concerned about what the effect of the prior is on the analysis (Evans, 2015). Methods for checking the prior in previous sense are developed in Evans and Moshonov (2006). See also Nott, Xueou, Evans, and Engler (2016)  and Nott, Seah, AL-Labadi, Evans, Ng and Englert (2019). The basic method for checking the prior involves computing the probability
\begin{eqnarray}
M_T \left(m_T(t) \le m_T(T(x))\right), \label{eq1}
\end{eqnarray}
where $T$ is a minimal sufficient statistic of the model and $M_T$ is
the prior predictive probability measure of $T$ with density $m_T$. The value of
(\ref{eq1}) simply serves to  locate the observed value $T(x)$ in its prior distribution.
If (\ref{eq1}) is small, then $T(x)$ lies in a region of low prior
probability, such as a tail or anti-mode, which indicates a conflict.
The consistency of this check follows from Evans and Jang (2011) where it is
proven that, under quite general conditions, (\ref{eq1}) converges to
\begin{eqnarray}
\Pi_T \left(\pi_0(\theta) \le \pi_0(\theta_{\text{true}})\right), \label{eq2}
\end{eqnarray}  as the amount of data increases, where $\theta_{\text{true}}$
is the true value of the parameter. If (\ref{eq2}) is small,
then $\theta_{\text{true}}$ lies in a region of low prior probability which implies
that the prior is not appropriate.


\section{A Bayesian Alternative to the One-Sample $z-$Test}

\subsection{The Approach} Let $x=(x_1,\ldots,x_n)$ be an independent random sample from $N(\mu, \sigma^2)$, where
$\sigma^2$ is known. The goal is to test the hypothesis $H_0:\mu = \mu_1$, where $\mu_1$
is a given real number. The approach here is  Bayesian. First we construct
a prior $\pi(\mu)$ on $\mu$. Let $\pi(\mu)$ be $N(\mu_0, \lambda^2_0 \sigma^2)$, where $\mu_0$ and
$\lambda_0^2$ are known hyperparameters and selected through the elicitation algorithms covered in Section 3.2. Thus, the posterior distribution of $\mu$ given $x_1,\ldots,x_n$ is
$\pi(\mu|x_1,\ldots,x_n)=N(\mu_x, \sigma_x^2)$, where
\begin{eqnarray}
\mu_x = \dfrac{n\lambda_0^2}{n \lambda_0^2 + 1}\bar{x} +
\dfrac{1}{n\lambda_0^2 + 1}\mu_0 \text{ and }
\sigma_x^2 = \dfrac{\lambda_0^2}{n\lambda_0^2 + 1}. \label{post}
\end{eqnarray}	
To proceed for the test using the relative belief ratio, there are two possible approaches. The first one is based on a direct computation of the relative belief ratio $RB(\mu_1|x)$ and its strength. This approach has been initiated in Baskurt and Evans (2013) with $\sigma^2=1$ and $\mu_1=0$ when discussing the Jeffrey-Lindely paradox. To find $RB(\mu|x)$, notice that
\begin{eqnarray}
\nonumber RB(\mu|x)=\frac{\pi(\mu|T(x))}{\pi(\mu)}=\frac{\pi(\mu)f(T(x))/m_T(T(x))}{\pi(\mu)}=\frac{f(T(x))}{m_T(T(x))}.
\end{eqnarray}
The minimal sufficient statistics for $\mu$ is $T(x)=\bar{x} \sim N(\,\mu,\sigma^2/n)$. Since  $T(x)=\bar{x}=(\bar{x}-\mu)+\mu$, where $\bar{x}-\mu \sim N(0, \sigma^2/n)$ independent of  $\mu \sim N(\mu_0,\lambda^2_0 \sigma^2)$, it follows the prior predictive distribution of $T(x)$ is $N(\mu_0, \lambda^2_0 \sigma^2_0+\sigma^2/n)$. That is,
\begin{eqnarray}
\nonumber m_T(T(x))=\sqrt{\frac{n}{2\pi\sigma^2(1+n\lambda_0^2)}}\exp\left(-\frac{n}{2}\frac{(\bar{x}-\mu_0)^2}{\sigma^2(n\lambda_0^2+1)}\right).
\end{eqnarray}
Thus,
\begin{eqnarray}
 RB(\mu|x)=\sqrt{1+n\lambda_0^2}\exp\left(-\frac{n}{2\sigma^2}\left[(\bar{x}-\mu)^2-\frac{(\bar{x}-\mu_0)^2}{(n\lambda_0^2+1)}\right]\right). \label{RB_Case1}
\end{eqnarray}
For the strength, we have $\Pi\left(RB(\mu|x) \le RB(\mu_1|x)|x\right)=$
\begin{eqnarray}
\nonumber &&\Pi\Bigg(\exp\left(-\frac{n}{2\sigma^2}\left[(\bar{x}-\mu)^2-\frac{(\bar{x}-\mu_0)^2}{(n\lambda_0^2+1)}\right]\right)\\
 \nonumber &&\le \exp\left(-\frac{n}{2\sigma^2}\left[(\bar{x}-\mu_1)^2-\frac{(\bar{x}-\mu_0)^2}{(n\lambda_0^2+1)}\right]\right)\Bigg|x\Bigg) \\
\nonumber &=&  \Pi\Bigg((\mu-\bar{x})^2 \ge (\bar{x}-\mu_1)^2\Bigg|x\Bigg)\\
\nonumber &=&  \Pi\bigg(|\mu-\bar{x}| \ge |\bar{x}-\mu_1|\Bigg|x\bigg)\\
\nonumber &=&\Pi\bigg(\mu \ge \bar{x}+|\bar{x}-\mu_1|\big|x\bigg) +\Pi\bigg(\mu \ge \bar{x}-|\bar{x}-\mu_1|\Bigg|x\bigg)\\
\nonumber &=&1-\Phi\bigg(\frac{\bar{x}+|\bar{x}-\mu_1|-\mu_x}{\sigma_x}\bigg) \\
\nonumber          &&+\Phi\bigg(\frac{\bar{x}-|\bar{x}-\mu_1|-\mu_x}{\sigma_x}\bigg),
\end{eqnarray}
where $\mu_x$ and $\sigma_x$ are defined in (\ref{post}). After minor simplification we have,
\begin{eqnarray}
\nonumber \Pi\left(RB(\mu|x) \le RB(\mu_1|x)|x\right)&=&1-\Phi\Bigg(\left(\frac{1}{\sigma^2}+\frac{1}{n \lambda_0^2\sigma^2}\right)^{1/2}\left(\sqrt{n}\left|\bar{x}-\mu_1\right|\right)\\
\nonumber                                                     &&+\frac{\sqrt{n}\bar{x}}{n\lambda_0^2+1}-\frac{\sqrt{n}\mu_0}{n\lambda_0^2+1}\Bigg)\\
\nonumber                                                           && +\Phi\Bigg(\left(\frac{1}{\sigma^2}+\frac{1}{n \lambda_0^2\sigma^2}\right)^{1/2}\left(-\sqrt{n}\left|\bar{x}-\mu_1\right|\right)\\
                                                    &&+\frac{\sqrt{n}\bar{x}}{n\lambda_0^2+1}-\frac{\sqrt{n}\mu_0}{n\lambda_0^2+1}\Bigg). \label{stength_dir}
\end{eqnarray}
Similar to the conclusion in Baskurt and Evans (2013),  as $\lambda_0^2 \to \infty$ in (\ref{stength_dir}), $\Pi\left(RB(\mu|x) \le RB(\mu_1|x)|x\right) \to 2\left(1-\Phi(\sqrt{n}|\bar{x}-\mu_1|/\sigma)\right)$, which converges in distribution to $2\left(1-\Phi(|z|)\right)$ when $\mu=\mu_1$, by the central limit theorem and the continuous mapping theorem, where $z$ is the standard normal random variable. Hence, when $\mu=\mu_1$ (i.e. $H_0$ is not rejected), the strength has an asymptotically uniform distribution on $(0,1)$. On the other hand, we have $\Pi\left(RB(\mu|x) \le RB(\mu_1|x)|x\right)$ converges to 0 almost surely (a.s.) when $\mu \neq \mu_1$, since $\bar{x} \to \mu$ almost surely.

As for the second approach, we compute the  KL distance between the  hypothesized distribution and the prior/posterior distributions. The change of the distance from a priori to a posteriori is compared through the relative belief ratio. Then,  we give a brief summary about the KL distance. In general, the KL distance (sometimes called the \emph{entropy distance}) between two  continuous cumulative distribution functions (cdf's) $P$ and $Q$ with corresponding probability density functions (pdf's)  $p$ and $q$ (with respect to Lebesgue measure) is defined by
\begin{eqnarray*}
d(P,Q)=\int p(x)\log\left(\frac{p(x)}{q(x)}\right) dx.
\end{eqnarray*}
 It is well-known that $d_{KL}(P,Q) \ge 0$ and the equality holds if and only if $p=q$.  However, it is not symmetric and does not satisfy the triangle inequality (Cover and Thomas, 1991).
In particular, the KL divergence between the two normal distributions $P=N(\mu_1, \sigma_1^2)$ and
$Q=N(\mu_2, \sigma_2^2)$ is given by (Duchi, 2014)
\begin{eqnarray}\label{KL1}
d(P, Q) = \log{(\tfrac{\sigma_1}{\sigma_2})} + \dfrac{1}{2\sigma_2^2}
\left[\sigma_1^2+(\mu_1 - \mu_2)^2\right] - \dfrac{1}{2}.
\end{eqnarray}
Set $P^{\text{prior}} = N(\mu, \sigma^2)$ and $Q = N(\mu_1, \sigma^2)$. It follows that from (\ref{KL1}) that
\begin{eqnarray}\label{prior1}
d(P^{\text{prior}}, Q) = \dfrac{(\mu - \mu_1)^2}{2\sigma^2}.
\end{eqnarray}
Also the KL divergence between $P^{\text{post}} = N(\mu_x, \sigma^2)$ and Q is
\begin{eqnarray}\label{post1}
d(P^{\text{post}},Q) = \dfrac{(\mu_x - \mu_1)^2}{2\sigma^2}.
\end{eqnarray}
Note that, as $n\rightarrow\infty$, by the strong law of large numbers, $\mu_x \overset{a.s.}\to \mu_{\text{true}}$, where $\mu_{\text{true}}$ is the true value of $\mu$. Thus, by (\ref{post1}), if $H_0$ is true,
we have $d(P^{\text{post}},Q) \overset{a.s.}\to 0$. On the other hand, if $H_0$ is  not true, then
\begin{eqnarray}\label{postnottrue}
	d(P^{\text{Post}}, Q) \overset{a.s.} \to c > 0.
\end{eqnarray}
What follows is that, if $H_0$ is true, then that distribution of $d(P^{\text{Post}},Q)$ should be more
concentrated about $0$ than $d(P^{\text{Prior}}, Q)$. So, the proposed test includes a comparison of the concentrations of the
prior and posterior distributions of the KL divergence via a relative belief ratio
based on  the interpretation as discussed in$\ $Section 2.

\subsection {Elicitation  of the Prior}
The success of methodology is influenced significantly by  the choice of the hyperparameters $\mu_0$ and $\lambda_0$.  Inappropriate values of the hyperparameters can lead to a failure in computing $d$. To elicit proper values of the hyperparameters, we consider the method developed in Evans and Tomal (2018). Suppose that it is known with virtual certainty, based on the knowledge of the basic measurement being taken, that $\mu$ will lie in the interval $(a, b)$ for specified values $a \le b$. Here, virtual certainty is interpreted as $P(a\le \mu \le b) \ge \gamma$,  where $\gamma$ is a large probability like 0.999. If $\mu_0=(a+b)/2$, then after some simple algebra,  $\lambda_0=(b-a)/(2\sigma \Phi^{-1}((1+0.999)/2))$.

\subsection {Checking for Prior-Data Conflict}

As pointed in Section 3.1, the minimal sufficient statistics for $\mu$ is $T(x)=\bar{x}$ with the prior predictive distribution of $T(x)$ is $N(\mu_0, \lambda^2_0 \sigma^2_0+\sigma^2/n)$. Thus,
 \begin{eqnarray}\label{conflict1}
	M_T\left(m_T(t) \le m_T(\bar{x})\right)=2\left(1-\Phi\left(\left|\bar{x}-\mu_0\right|/\left(\lambda^2_0\sigma^2+\sigma^2/n\right)^{0.5}\right)\right),
\end{eqnarray}
where $M_T$ is defined as in (\ref{eq1}). Recall that, if (\ref{conflict1}) is small, then this indicates a prior-data conflict and no prior-data conflict otherwise. It is true that prior-data conflict can be avoided by increasing $\lambda_0$ (i.e. making the prior diffuse),  however, as pointed in Evans (2018), this is not an appropriate approach as it will induce bias into the analysis. Thus, by (\ref{conflict1}), when $\bar{x}_0$ lies in the tail of its prior distribution, we have a prior-data conflict. Note that, as $n \to \infty,$ $(\ref{conflict1})\overset{a.s.} \to 2\left(1-\Phi\left(\left|\mu_{\text{true}}-\mu_0\right|/\left(\lambda_0\sigma\right)\right)\right)$.

\subsection {Checking for Bias}
The bias against the hypothesis $H_0: \mu=\mu_1$ is measured by computing (\ref{bias_against}) with $\psi_0=\mu_1$ and $RB(\mu_1|x)$ as  in (\ref{RB_Case1}).  Note that, since the prior is centered at $\mu_1$, there is never a strong bias against $H_0$.  On the other hand, the bias in favor of the hypothesis $H_0: \mu=\mu_1$ is measured by computing (\ref{bias_favor}) with $\psi_0=\mu_1$ and $RB(\mu_1|x)$ as defined in (\ref{RB_Case1}). The interpretation of the bias was covered in Section 2.

\subsection {The Algorithm}
The approach will involve a comparison between the concentrations of the prior and posterior distribution of the KL
divergence via a relative belief ratio, with the interpretation as discussed in Section 2. Since explicit forms of the densities of the distance are not available, the relative belief
ratios need to be estimated via simulation. The following summarizes a
computational algorithm for testing $H_0$.\\
		
\underline{Algorithm A (New $z-$Test)}
\begin{enumerate} [label=(\roman*)]
\item  Elicit the hyperparameters $\mu_0$ and $\lambda_0$ as described in Section 3.2.
\item Generate $\mu$ from $N(\mu_0, \lambda^2_0\sigma_0^2)$.
\item Compute the KL distance between $N(\mu, \sigma^2)$ and $Q=N(\mu_1, \sigma^2)$ as described in (\ref{prior1}). Denote this distance by $D$.
\item Repeat steps (ii) and (iii) to obtain a sample of $r_1$ values of $D$.
\item Generate $\mu$ from $N(\mu_x, \sigma^2_x)$, where $\mu_x$ and  $\sigma^2_x$ are defined in (\ref{post}).
\item Compute the KL distance between $N(\mu, \sigma^2)$ and $Q=N(\mu_1, \sigma^2)$ as described in (\ref{post1}). Denote this distance by $D_x$.
\item Repeat steps (v) and (vi) to obtain a sample of $r_2$ values of $D_x$.
\item Compute the relative belief ratio and the strength as follows:
\begin{enumerate}

\item  Closed forms of $D$ and $D_x$ are not available. Thus, the relative brief ration and the strength need to be estimated via approximation.  Let $M$ be a positive number. Let $\hat{F}_{D}$ denote the
empirical cdf of $D$ based on the prior sample in (3) and for $i=0,\ldots,M,$
let $\hat{d}_{i/M}$ be the estimate of $d_{i/M},$ the $(i/M)$-the prior
quantile of $D.$ Here $\hat{d}_{0}=0$, and $\hat{d}_{1}$ is the largest value
of $d$. Let $\hat{F}_{D}(\cdot\,|\,x)$ denote the empirical cdf of $D$ based
on the posterior sample in (vi). For $d\in\lbrack\hat{d}_{i/M},\hat
{d}_{(i+1)/M})$, estimate $RB_{D}(d\,|\,x)={\pi_D(d|x)}/{\pi_D(d)}$ by
\begin{equation}
\widehat{RB}_{D}(d\,|\,x)=M\{\hat{F}_{D}(\hat{d}_{(i+1)/M}\,|\,x)-\hat{F}%
_{D}(\hat{d}_{i/M}\,|\,x)\}, \label{rbest}%
\end{equation}
the ratio of the estimates of the posterior and prior contents of $[\hat
{d}_{i/M},\hat{d}_{(i+1)/M}).$ Thus,   we estimate $RB_{D}(0\,|\,x)={\pi_D(0|x)}/{\pi_D(0)}$
 by
$\widehat{RB}_{D}(0\,|\,x)=$ $M\widehat{F}_{D}(\hat{d}_{p_{0}}\,|\,x)$ where
$p_{0}=i_{0}/M$ and $i_{0}$ are chosen so that $i_{0}/M$ is not too small
(typically $i_{0}/M\approx0.05)$.\textbf{\smallskip}

\item Estimate the strength $DP_{D}(RB_{D}(d\,|\,x)\leq RB_{D}%
(0\,|\,x)\,|\,x)$ by the finite sum
\begin{equation}
\sum_{\{i\geq i_{0}:\widehat{RB}_{D}(\hat{d}_{i/M}\,|\,x)\leq\widehat{RB}%
_{D}(0\,|\,x)\}}(\hat{F}_{D}(\hat{d}_{(i+1)/M}\,|\,x)-\hat{F}_{D}(\hat
{d}_{i/M}\,|\,x)). \label{strest}%
\end{equation}

\noindent For fixed $M,$ as $r_{1}\rightarrow\infty,r_{2}\rightarrow\infty,$
then $\hat{d}_{i/M}$ converges almost surely to $d_{i/M}$ and (\ref{rbest})
and (\ref{strest}) converge almost surely to $RB_{D}(d\,|\,x)$ and
$DP_{D}(RB_{D}(d\,|\,x)\leq RB_{D}(0\,|\,x)\,|\,x)$, respectively.

\end{enumerate}

\end{enumerate}

The following proposition establishes the consistency of the approach
as the sample size increases. So, the procedure performs correctly as the sample size increases when
$\mathcal{H}_{0}$ is true. The proof follows immediately from Evans (2015), Section 4.7.1. See also AL-Labadi and Evans (2018) for a similar result.

\begin{proposition}
\label{cvm6}Consider the discretization $\{[0,d_{i_{0}/M}),[d_{i_{0}%
/M},d_{(i_{0}+1)/M}),\ldots,$\newline$[d_{(M-1)/M},\infty)\}$. As
$n \rightarrow\infty,$ (i) if $\mathcal{H}_{0}$ is true, then
\begin{align*}
&  RB_{D}([0,d_{i_{0}/M})\,|\,x)\overset{a.s.}{\rightarrow}1/DP_{D}%
([0,d_{i_{0}/M})),\\
&  RB_{D}([d_{i/M},d_{(i+1)/M})\,|\,x)\overset{a.s.}{\rightarrow}0\text{
whenever }i\geq i_{0},\\
&  DP_{D}(RB_{D}(d\,|\,x)\leq RB_{D}(0\,|\,x)\,|\,x)\overset{a.s.}%
{\rightarrow}1,
\end{align*}
and (ii) if $\mathcal{H}_{0}$ is false and $d_{CvM}(P, Q)\geq d_{i_{0}/M}$, then $RB_{D}([0,d_{i_{0}/M})\,|\,x)\overset
{a.s.}{\rightarrow}0$ and $DP_{D}(RB_{D}(d\,|\,x)\leq RB_{D}%
(0\,|\,x)\,|\,x)\overset{a.s.}{\rightarrow}0.$
\end{proposition}

%

\section{A Bayesian Alternative to the One-Sample t-Test}
\subsection{The Approach}
In this section, we assume that $x_1, ..., x_n$ is an independent random sample from $\mathcal{N}(\mu, \sigma^2)$, where $\sigma^2$ is unknown. The goal is to test $H_0: \mu = \mu_1$, where $\mu_1$ is a given real number. The first step in the approach is to construct priors on $\mu$ and $\sigma^2$. We will consider the following hierarchical but conjugate prior (Evans 2015, p.171):
\begin{align}
    \frac{1}{\sigma^2} \sim \text{gamma}_{rate}(\alpha_0, \beta_0) \label{prior2a}
\\
    \mu | \sigma^2 \sim \mathcal{N}(\mu_0, \lambda_0^2\sigma^2), \label{prior2b}
\end{align}
where $\mu_0$, $\lambda_0$ and $(\alpha_0, \beta_0)$ are hyperparameters to be specified via elicitation as it will be described in Section 4.2.
The posterior distribution of $(\mu, \sigma^2)$ is given by:
\begin{align}
   \frac{1}{\sigma^2} \big| x_1, ..., x_n \sim \text{gamma}_{rate} (\alpha_0 + \frac{n}{2}, \beta_x), \label{post2a}
 \\
\mu | \sigma^2, x_1, ..., x_n \sim \mathcal{N}\left(\mu_x, (n + \frac{1}{\lambda_0^2})^2 \sigma^2\right) \label{post2b}
\end{align}
where
\begin{equation}
\mu_x = \left(n + \frac{1}{\lambda_0}\right)^{-1} \left(\frac{\mu_0}{\lambda_0^2} + n\bar{x}\right) \text{ and } \beta_x = \beta_0 + (n-1)\frac{S^2}{2} + \frac{n(\bar{x} - \mu_0)^2}{2(n\lambda_0^2 + 1)} \label{beta_x}
\end{equation}
with $S^2 = \frac{1}{n-1}\sum_{i=1}^n (x_i - \bar{x})^2.$  To find $RB(\mu|x)$, notice that the minimal sufficient statistic for $T(x)=(\mu, \sigma^2)$ is $(\bar{x}, s^2)$ with $\bar{x} \sim \mathcal{N}\left(\mu, \sigma^2/n\right)$ independent of $s^2 \sim \sigma^2(n-1)^{-1}\chi^2_{n-1}$. The joint prior predictive of $T(x) = (\bar{x}, s^2)$ is given by (Evan, 2015):
\begin{align}
    m_T(T(x)) = \frac{\Gamma(\frac{n}{2} + \alpha_0)}{\Gamma(\alpha_0)}(n+\frac{1}{\lambda_0^2})^{-\frac{1}{2}}
    \frac{(2\pi)^{-\frac{n}{2}}\beta_0^{\alpha_0}}{\lambda_0}(\beta_{x})^{-\frac{n}{2} - \alpha_0}, \label{eq23}
\end{align}
where $\beta_x$ is defined in (\ref{beta_x}). On the other hand, it can be shown that
\begin{align*}
    m_T(T(x)|\mu) = \frac{\Gamma(\frac{n}{2} + \alpha_0)}{\Gamma(\alpha_0)}
    (2\pi)^{-\frac{n}{2}}\beta_0^{\alpha_0} \left(\beta_0 + \frac{n-1}{2}s^2 + \frac{n}{2}(\bar{x} - \mu)^2\right)^{-\frac{n}{2} - \alpha_0}.
\end{align*}
Thus,
\begin{align}
 \nonumber    RB(\mu|x) &= \frac{m_T(T(x)|\mu)}{m_T(T(x))}
   \\
 \nonumber    &= \left( n + \frac{1}{\lambda_0^2} \right) ^ {\frac{1}{2}}
    \frac{\left[
        \beta_0 + \frac{n-1}{2}s^2 + \frac{n}{2}(\bar{x} - \mu)^2
    \right]^{-\frac{n}{2} - \alpha_0} } {(\beta_{x})^{-\frac{n}{2} - \alpha_0}}
    \\
&= \left( n + \frac{1}{\lambda_0^2} \right) ^ {\frac{1}{2}}
    \left[
    \frac{
        \beta_0 + \frac{n-1}{2}s^2 + \frac{n}{2}(\bar{x} - \mu)^2
    }
    {
        \beta_0 + \frac{n-1}{2}s^2 + \frac{n}{2}\frac{(\bar{x} - \mu)^2}{n\lambda_0^2 + 1}
    }
    \right]^{-\frac{n}{2} - \alpha_0}. \label{eq25}
\end{align}
For the strength we have, $\Pi\left(RB(\mu|x) \le RB(\mu_1|x)|x\right)=$
\begin{eqnarray}
\nonumber &=&\Pi\left(
    \frac{\beta_0 + \frac{n-1}{2}s^2 + \frac{n}{2}(\bar{x} - \mu)^2} { \beta_0 + \frac{n-1}{2}s^2 + \frac{n}{2}\frac{(\bar{x} - \mu)^2}{n\lambda_0^2 + 1}
    }
\le
    \frac{\beta_0 + \frac{n-1}{2}s^2 + \frac{n}{2}(\bar{x} - \mu_1)^2 } {\beta_0 + \frac{n-1}{2}s^2 + \frac{n}{2}\frac{(\bar{x} - \mu)^2}{n\lambda_0^2 + 1}} \bigg|x\right),\\
\end{eqnarray}
where $\mu_x$ and $\sigma_x$ are defined in (\ref{post2a}) and (\ref{post2b}), respectively. After some algebra, we reach the conclusion that $\Pi\left(RB(\mu|x) \le RB(\mu_1|x)|x\right)$ coincides with (\ref{stength_dir}), but here $\sigma^2$ is random as defined in (\ref{prior2a}).

As for  the  KL approach, we compute $d(p^{prior}, Q)$ and $d(p^{post}, Q)$ as  given respectively in (\ref{prior1}) and (\ref{post1}). The approach makes a comparison between the concentrations of the prior and posterior distributions of the KL divergence via the relative belief ratio.

\subsection{Elicitation of the prior}
To elicit the prior, we consider the approach developed by Evan (2015, p.171). Suppose that it is known with virtual certainty (probability = 0.999) that $\mu \in (a, b)$ for specified values $a \le b$. This is chosen to be as short as possible, based on the knowledge of the basic measurements being taken and without being unrealistic. We set $\mu_0 = (a + b)/2$ (i.e, mid-point). With this choice, one hyper-parameter has been specified. It follows that
\begin{align*}
    P(\mu \in (a, b)) \geq 0.999&\implies P(a < \mu < b) \geq 0.999
    \\
    &\implies P\left(\frac{a - \mu_0}{\lambda_0 \sigma_0} < Z < \frac{b - \mu_0}{\lambda_0 \sigma_0}\right) \geq 0.999
    \\
    &\implies \Phi\left(\frac{b - \mu_0}{\lambda_0 \sigma_0}\right) -
    \Phi\left(\frac{a - \mu_0}{\lambda_0 \sigma_0}\right) \geq 0.999
    \\
   & \implies 2\Phi\left(\frac{b - a}{2\lambda_0 \sigma_0}\right) - 1 \geq 0.999.
\end{align*}
This implies that
\begin{align}
  \nonumber  &\Phi\left(\frac{b - a}{2\lambda_0 \sigma_0}\right) \geq \frac{1.99}{2}=0.9995
    \\
    \nonumber   &\implies \frac{b - a}{2\lambda_0 \sigma_0} \geq
    \Phi^{-1}(0.9995)
    \\
   \nonumber    &\implies \sigma \leq \frac{b - a}{2\lambda_0 \Phi^{-1}(0.9995)}
    \\
    &\implies \sigma^2 \leq \left(\frac{b - a}{2}\right)^2 \left[\Phi^{-1}(0.9995)\right]^{-2} \lambda_0^{-2}. \label{eq18}
\end{align}
An interval that contains virtually all the actual data measurements is given by $\mu \pm \sigma \Phi^{-1} (0.9995)$. Since this interval cannot be unrealistically too short or too long, we let $s_1$ and $s_2$ be the upper and lower bounds on the half-length of the interval so that
\begin{align*}
    s_1 \leq \sigma \Phi^{-1} (0.9995) \leq s_2.
\end{align*}
That is,
\begin{align}
 \frac{s_1}{\Phi^{-1} (0.9995)} \leq \sigma  \leq \frac{s_2}{\Phi^{-1} (0.9995)} \label{eq19}
\end{align}
Now, from (\ref{eq18}) and (\ref{eq19}), we have:
\begin{align*}
    \left(\frac{b - a}{2}\right)^2 \left[\Phi^{-1}(0.9995)\right]^{-2} \lambda_0^{-2}
    = \left(\frac{s_2}{\Phi^{-1} (0.9995)}\right)^2
    \\\\
    \implies \lambda_0^2 = \left(\frac{b - a}{2}\right)^2  s_2^{-2},
\end{align*}
which determine the conditional prior for $\mu$. Note that $\lambda_0$ can be made bigger by choosing a bigger value of $b - a$.

Lastly, to obtain relevant values of $\alpha_0$ and $\beta_0$, let $G(\alpha_0, \beta_0, x)$ denotes the CDF of $\text{gamma}_{rate}(\alpha_0, \beta_0)$ distribution.  From(\ref{eq19}),
\begin{align}
    s_2^2\left[\Phi^{-1}(0.9995)\right] \leq
    \frac{1}{\sigma^2} \leq
    s_1^2\left[\Phi^{-1}(0.9995)\right]. \label{eq20}
\end{align}
Now, suppose we want to determine the lower and upper bounds in (\ref{eq20}), so that this interval contains $1/\sigma^2$ with virtual certainty. Thus,
\begin{align}
    G^{-1}(\alpha_0, \beta_0, 0.9995) = s_1^{-2}\left[\Phi^{-1}(0.9995)\right]^2 \label{eq21}
    \\
    G^{-1}(\alpha_0, \beta_0, 0.0005) = s_2^{-2}\left[\Phi^{-1}(0.9995)\right]^2. \label{eq22}
\end{align}

Then we numerically solve (\ref{eq21}) and (\ref{eq22})  for $(\alpha_0, \beta_0)$.


\subsection{Checking for Prior-data Conflict}
To assess whether $(\bar{x}_0, s_0^2)$ is a reasonable value, we compute:
\begin{align}
    M_T\left(m_T(\bar{x}, s^2) \leq m_T(\bar{x}_0, s_0^2)\right),
\label{eq24}
\end{align}
where $T$, $M_T$ and $m_T$ are as defined in Section 4.1. Clearly, computing (\ref{eq24}) should be done by simulation. Thus, for specified values of $\mu_0, \lambda_0^2, (\alpha_0, \beta_0)$, we generate $(\mu, \sigma^2)$ as given in (\ref{prior2a}) and (\ref{prior2b}). Then generate $(\bar{x}, s^2)$ from the joint distribution given $(\mu, \sigma^2)$ and evaluate $m_T(\bar{x}, s^2)$ using (\ref{eq23}). Repeating this many times and recording the proportion of values of $m_T(\bar{x}, s^2)$ that are less than or equal to $m_T(\bar{x}_0, s_0^2)$ gives a Monte Carlo estimate of (\ref{eq24}).


\subsection{Checking for Bias}
As in Section 3.4, the bias against the hypothesis $H_0: \mu = \mu_1$ is measured by computing (\ref{bias_against})  with $\psi_0=\mu_1$ and   $RB(\mu|x)$  as given in (\ref{eq25}). On the other hand, the bias in favor of the hypothesis $H_0: \mu = \mu_1$ is measured by computing (\ref{bias_favor})  with $\psi_0=\mu_1$ and $RB(\mu_1|x)$ as defined in (\ref{eq25}). The interpretation of the bias was given in Section 2.

\subsection {The Algorithm}
The following algorithm outlines the KL approach described in Section 4.
		
\underline{Algorithm B (New $t-$Test)}
\begin{enumerate} [label=(\roman*)]
\item  Elicit the hyperparameters $\mu_0$, $\lambda_0$ and $(\alpha_0,\beta_0)$ as described in Section 4.2.
\item Generate $\mu$ and $\sigma^2$ as described in (\ref{prior2a}) and (\ref{prior2b}).
\item Compute the KL distance between $P^{\text{Prior}}=N(\mu, \sigma^2)$ and $Q=N(\mu_1, \sigma^2)$ as described in (\ref{prior1}). Denote this distance by $D$.
\item Repeat steps (ii) and (iii) to obtain a sample of $r_1$ values of $D$.
\item  Generate $\mu_x$ and $\sigma^2_x$ from (\ref{post2a}) and (\ref{post2b}), respectively.
\item Compute the KL distance between $P^{\text{Post}}=N(\mu, \sigma^2)$ and $Q=N(\mu_1, \sigma^2)$ as described in (\ref{post1}). Denote this distance by $D_x$.
\item Repeat steps (v) and (vi) to obtain a sample of $r_2$ values of $D_x$.
\item Compute the relative belief ratio and the strength described in Algorithm A.

\end{enumerate}

Note that, like Proposition \ref{cvm6}, the approach in this  case (i.e., when $\sigma^2$ is unknown) is also consistent as the sample size increases.

\section{Examples}
In this section, we consider three examples. The first one deals with a study on dental anxiety in adults, where the goal  is to gauge the fear of adults of going to a dentist (McClave and Sincich, 2017, p. 398). For this, a random sample of  15  adults completed the Modified Dental Anxiety Scale questionnaire, where scores range from zero (no anxiety) to 25 (extreme anxiety). The sample mean score was 10.7 and the sample standard deviation was 3.6. We want to determine whether the mean Dental Anxiety  Scale score for the population differs from 11. To construct the prior, we implement the elicitation algorithm described in Section 4.1 with $a=0$, $b=25$,  $s_1=2$, $s_2=15$  and $\gamma=0.999$. Consequently, we have $\mu_0=12.5$, $\lambda_0= 0.83$, $\alpha_0=1.29$ and $\beta_0=12.36$. To check if there is a prior-data conflict, (\ref{eq24}) is computed to be 0.46, and that implies an indication of no prior-data conflict.   The bias is also assessed by computing (\ref{bias_against}) with $\psi_0=11$.  In this case the bias against the null hypothesis is  0.5136.  On the other hand, the bias in favor of the null hypothesis  is measured by computing (\ref{bias_favor})  with $\psi_0=11 \pm 0.5$, which gives 0.5192  (for $\psi_0=11.5$) and 0.5208 (for $\psi_0=10.5$). This shows equal  bias either for or against the null hypothesis for this choice of prior. The value of the relative believe ratio test (distance and direct) with strength, the test of Rouder et.\ al.\ (2009) and the standard $t$-test are summarized in Table 1. It follows from the table that the null hypothesis is accepted by the three  Bayesian  tests, while it is not rejected by the $t$-test.

\begin{table}[htp] \centering
\begin{tabular}
[c]{|l|l|l|}\hline
\multicolumn{1}{|c|}{Test}  &  \multicolumn{1}{|c|} {Values} & \multicolumn{1}{|c|}{Decision}\\ \hline
Distance: RB (Strength) & 4.7800(1)& Accept the null hypothesis\\ \hline
Direct: RB (Strength) &  4.8466(0.4341)&Accept the null hypothesis\\ \hline
Rouder et. al. (2009): BF& 0.2747& Accept the null hypothesis\\  \hline
$t$-test: p-value& 0.7517 & Fail to reject the null hypothesis \\\hline
\end{tabular}
\caption{The tests results about the dental anxiety example.}%
\end{table}%

The second example considers an application about the age at which children start walking (Mann, 2016). A psychologist claims that the mean age at which children start walking is 12.5 months. To test this claim, she took a random sample of 18 children and found that the mean age  at which they started walking was 12.9 with a standard deviation of 0.80 month. As in the previous example, the prior is constructed by setting $a= 8$, $b=24$,  $s_1=4$, $s_2=10$  and $\gamma=0.999$ the algorithm described in Section 4.1. It follows that $\mu_0=16$, $\lambda_0= 0.8$, $\alpha_0= 4.01$ and $\beta_0= 329.78$. With this prior, it is found that (\ref{eq24})=1, which is a clear indication of no prior-data conflict.   The bias is also assessed by computing (\ref{bias_against})  with $\psi_0=12.5$. In this case the bias against the null hypothesis is  0.4991. Moreover, the bias in favor of the null hypothesis  is measured by computing (\ref{bias_favor})  with $\psi_0=12.5 \pm 0.5$, which gives 0.5169  (for $\psi_0=13$) and  0.5140 (for $\psi_0=12.0$). This shows equal  bias either for or against the null hypothesis for this choice of prior. The results are reported in Table 2. Thus, the tests based on the relative belief ratio accept the null hypotheses  while the other two tests reject the null hypothesis.

\begin{table}[htp] \centering
\begin{tabular}
[c]{|l|l|l|}\hline
\multicolumn{1}{|c|}{Test}  &  \multicolumn{1}{|c|} {Values} & \multicolumn{1}{|c|}{Decision}\\ \hline
Distance: RB (Strength) &  6.7480(1)& Accept the null hypothesis\\ \hline
Direct: RB (Strength) &  6.5420(0.4479)&Accept the null hypothesis\\ \hline
Rouder et. al. (2009): BF& 1.467772& Reject the null hypothesis\\  \hline
$t$-test: p-value& .0489 & Reject the null hypothesis \\\hline
\end{tabular}
\caption{The tests results about the age at which children start walking.}%
\end{table}%

The last example deals with sugar production (Bluman, 2012, p. 457), where sugar is packed in 5-pound bags. An inspector suspects the bags may not contain 5 pounds. A sample of 50 bags produces a mean of 4.6 pounds and a standard deviation of 0.7 pound. To goal is to test if  bags do not contain 5 pounds as stated. In the algorithm given in Section 4.1, we set  $a= 4$, $b=6$,  $s_1=2$, $s_2=5$  and $\gamma=0.999$. We get $\mu_0=5$, $\lambda_0= 0.2 $, $\alpha_0= 4.0077$ and $\beta_0= 20.6106$. For this prior,  (\ref{eq24})=0.6262, which means no prior-data conflict is found.   The bias is  evaluated by computing (\ref{bias_against})  with $\psi_0=4.6$.  In this case, the bias against the null hypothesis is  0.4877. Additionally, the bias in favor of the null hypothesis  is measured by computing (\ref{bias_favor})  with $\psi_0=4.6 \pm 0.5$, which gives 0.4876  (for $\psi_0=5.1$) and  0.4864 (for $\psi_0=4.1$). This demonstrate equal bias either for or against the null hypothesis for this choice for prior. The results are reported in Table 3. From the table, we see that all the previous tests reject the null hypothesis.

\begin{table}[htp] \centering
\begin{tabular}
[c]{|l|l|l|}\hline
\multicolumn{1}{|c|}{Test}  &  \multicolumn{1}{|c|} {Values} & \multicolumn{1}{|c|}{Decision}\\ \hline
Distance: RB (Strength) &  0.4680(0.0442)& Reject the null hypothesis\\ \hline
Direct: RB (Strength) &  0.4355(0.0189)&Reject the null hypothesis\\ \hline
Rouder et. al. (2009): BF& 129.1731& Reject the null hypothesis\\  \hline
$t$-test: p-value& 0.0002 & Reject the null hypothesis \\\hline
\end{tabular}
\caption{The tests results about the sugar production.}%
\end{table}%

\section{Concluding Remarks}

A Bayesian approach to the standard one-sample $z$- and $t$- tests has been developed.   The prior has been created through an elicitation algorithm. Then the prior is evaluated for the existence of prior-data conflict and bias. The use of the approach has been illustrated through several examples, in which it shows excellent performance.

The approach can be extended in several directions. For instance, it can be used to test the difference between two population means. Also, it can be modified to be a Bayesian alternative to the Hotelling's  $T^2$ test for the multivariate mean.


\begin{thebibliography}{}


\bibitem{1} Al-Labadi, L. (2018).  The two-sample problem via relative belief ratio. https://arxiv.org/abs/1805.07238

\bibitem{2} Al-Labadi, L. and Evans, M. (2018). Prior based model checking. To appear in \emph{Canadian Journal of Statistics}.

\bibitem{4} Al-Labadi, L. and Evans, M. (2017). Optimal robustness results for relative belief inferences and the relationship
to prior-data conflict. \emph{Bayesian Analysis}, 12, 705--728.

\bibitem{5}  Al-Labadi, L., Baskurt , Z. and Evans, M. (2017). Goodness of fit for the logistic regression model using relative belief.   \emph{Journal of Statistical Distributions and Applications}. DOI 10.1186/s40488-017-0070-7.


\bibitem{3}  Al-Labadi, L., Baskurt , Z.  and Evans, M. (2018). Statistical reasoning: choosing and checking
the ingredients, inferences based on a measure of statistical evidence with some applications. \emph{Entropy}, 20, 289; doi:10.3390/e20040289.







\bibitem{6} Baskurt, Z.  and Evans, M. (2013). Hypothesis assessment and
inequalities for Bayes factors and relative belief ratios. \emph{Bayesian Analysis},
8, 3, 569-590.


\bibitem{7} Bluman, A. (2012). \emph{Elementary Statistics: A Step by Step Approach}, 8th edition. McGraw-Hill Education.

\bibitem {bib7}Cover, T. M. and  Thomas, J. A. (1991).  \emph{Elements of Information Theory}. Wiley.


\bibitem{8} Duchi, J. (2007). Derivations for Linear Algebra and Optimization. University of
California, Berkeley. URL: {\url{https://pdfs.semanticscholar.org/333f/05838e9b5a042b1b5b6113a404b74d22bedd.pdf}}

\bibitem{9} Evans, M. (2015). \emph{Measuring Statistical Evidence Using
Relative Belief}. Monographs on Statistics and Applied Probability 144, CRC
Press, Taylor \& Francis Group.

\bibitem{10} Evans, M. and Moshonov, H. (2006). Checking for prior-data conflict.
\emph{Bayesian Analysis}, 1,  893--914.

\bibitem{11} Evan, M. and Tomal, J. (2018). Measuring statistical evidence and multiple testing. \emph{FACET}, 3, 563--583.


\bibitem {12} G\"onen, M., Johnson, W. O.,  Lu, Y.  and  Westfall, P. H. (2005). The Bayesian two-sample t
test. \emph{The American Statistician}, 59, 252--25.

\bibitem{13} Kass, R. E. and Raftery,  A. E. (1995). Bayes factors. \emph{Journal of the American Statistical
Association}, 90, 773–-795.

 \bibitem{14} McClave, J. and Sincich, T. (2017). \emph{Statistics}, 13th ed. Person.

\bibitem {15} Nott, D. J., Xueou, W., Evans, M. and Englert, B. (2016). Checking for prior-data coflict using prior to posterior divergences. URL: {\url{https://arxiv.org/pdf/1611.00113.pdf}}.


\bibitem {16} Nott, D. J.,  Seah, M., AL-Labadi, L., Evans, M., Ng, H. K.  and Englert, B. (2019). Using prior expansion for prior-data cnflict checking.

\bibitem {17} Rouder, J., Speckman,  P. Sun,  D.,  Morey, R. and G. Iverson (2009). Bayesian t tests
for accepting and rejecting the null hypothesis. \emph{Psychonomic Bulletin \& Review}, 16, 225--23.

\bibitem{18} Rudin, W. (1974). \emph{Real And Complex Analysis}, 2nd ed. McGraw-Hill. 

\bibitem{19} Wang, M. and Lui, G. (2016). A simple two-sample Bayesian t-test for hypothesis testing.  \emph{The American Statistician}, 70, 195--201.


\bibitem{20} Wasserstein, R. L., and  Lazar, N. A. (2016). The ASA’s statement on p-values: Context, process, and purpose.
\emph{The American Statistician}, 70, 129--133. doi:10.1080/00031305.2016.1154108


\end{thebibliography}
\end{document}